\documentclass[titlepage,11pt]{article}
\usepackage{latexsym}
\usepackage{amsfonts}
\usepackage{amsmath}
\usepackage{graphics}
\usepackage{graphicx}

\newcommand\blackslug{\hbox{\hskip 1pt \vrule width 4pt height 8pt depth 1.5pt
        \hskip 1pt}}
\newcommand\bbox{\hfill \quad \blackslug \bigbreak}

\title{Tournament minors}
\author{Ilhee Kim and Paul Seymour\thanks{Supported by ONR grant N00014-10-1-0608 and NSF grant DMS-0901075.}\\
Princeton University, Princeton, NJ 08540}

\date{June 19, 2011; revised \today}

\newtheorem{thm}{}[section]

\newcommand{\Proof}{\noindent{\bf Proof.}\ \ }

\begin{document}
\maketitle

\begin{abstract}
We say a digraph $G$ is a {\em minor} of a digraph $H$
if $G$ can be obtained from a subdigraph of $H$ by repeatedly
contracting a strongly-connected subdigraph to a vertex.
Here, we show the class of all tournaments is a well-quasi-order under minor containment.
\end{abstract}

\section{Introduction}
The ``minor'' relation for graphs is well-established, but how it should be extended to digraphs is not clear.
In this paper we discuss one such extension. To obtain a minor, 
we contract strongly-connected subdigraphs rather than edges.
(A digraph $G$ is {\it strongly-connected} if $G$ is non-null and there exists a directed path from $u$ to $v$ for every $u,v \in V(G)$.)

In digraphs, contracting an edge may yield a directed cycle, even starting from an acyclic digraph, and this seems undesirable for a theory of excluded minors.
One way to avoid this is to only permit the contraction of certain edges.
For example, if an edge $uv$ is the only edge with tail $u$ or the only edge with head $v$, then contracting $uv$ does not 
yield a new directed cycle (see for instance \cite{DTree}).
Another way is to define a minor relation using minor mapping with certain conditions. 
See for instance \cite{Nowhere}.

A third way to extend minors of graphs to digraphs is as follows. For graphs, one can define contraction in terms of contracting edges, 
or in terms of contracting connected subgraphs, and it comes to the same thing. But for digraphs, the analogous two concepts are different. For digraphs, let us define contraction
in terms of contracting strongly-connected subgraphs.
More precisely, we say a digraph $H$ is a {\em minor} of a digraph $G$
if $H$ can be obtained from a subdigraph of $G$ by repeatedly
contracting a strongly-connected subdigraph to a vertex. (Note that we do not create ``new'' directed cycles after contracting a strongly-connected subdigraph.)
Equivalently, we can define this minor relation via minor mapping in the following way.
A digraph $H$ is a {\em minor} of a digraph $G$
if there exists a mapping $\phi$ defined on $V(H)$ such that:
\begin{itemize}
\item for every $v \in V(H)$, $\phi(v)$ is a non-null strongly-connected subdigraph of $G$,
\item if $u, v \in V(H)$ and $u \neq v$, then $\phi(u)$ and $\phi(v)$ are vertex-disjoint,
\item for every $u, v \in V(H)$ (not necessarily distinct),
if there are $k$ edges in $H$ with tail $u$ and head $v$,
then there are at least $k$ edges in $G$ with head in $V(\phi(u))$ and tail in $V(\phi(v))$,
and not contained in $E(\phi(x))$ for any $x\in V(H)$.
\end{itemize}
We call such a map $\phi$ a {\it minor mapping} from $H$ to $G$.

We first give some definitions.
Every digraph in this paper is finite.
We say a digraph $G$ is {\it simple} if it is loopless and
there is at most one edge $uv \in E(G)$ for every distinct $u, v \in V(G)$.
A simple digraph $G$ is called {\it semi-complete} if
either $uv \in E(G)$ or $vu \in E(G)$ for every distinct $u, v \in V(G)$.
A semi-complete digraph $G$ is called a {\it tournament} if
exactly one of $uv$ and $vu$ is an edge of $G$ for every distinct $u,v \in V(G)$.

An important property of minors for graphs is that they define a ``well-quasi-order''~\cite{GM20}.
A {\it quasi-order} $Q = (E(Q), \leq_Q)$ consists of a class $E(Q)$ and
a reflexive transitive relation $\leq_Q$ on $E(Q)$.
A quasi-order $Q$ is called a {\it well-quasi-order} or {\it wqo}
if for every infinite sequence $q_1, q_2, \ldots$ of elements of
$E(Q)$, there exist $j>i\ge 1$ such that $q_i \leq_Q q_j$.
Neil Robertson and the second author
proved that the class of all graphs is a wqo under the minor relation in
\cite{GM20}.
Unfortunately, the analogous statement is not true for directed minors.
For example, a small directed cycle is not a minor of a big directed cycle.
However, what if we consider some subclass, say the class of all tournaments?
(The subgraph relation does not define a wqo even for the class of all tournaments.
We leave it as an exercise for the reader to find a counterexample.)

In this paper, we prove that minor containment defines a wqo for the class of all semi-complete digraphs.
Therefore the same is true for the class of all tournaments.
We also give counter-examples for some other classes containing all semi-complete digraphs.

\begin{thm}\label{main1}
The class of all semi-complete digraphs is a wqo under minor containment.
\end{thm}

In \cite{WQO}, Maria Chudnovsky and the second author proved that
the class of all semi-complete digraphs is a wqo under
``immersion'', by using a digraph parameter called ``cut-width''.
Here, we prove the analogous statement for minors by using another parameter called {\it path-width}.
Path-width for undirected graphs was introduced in \cite{GM1}, and it has a natural
extension to digraphs, discussed for instance in~\cite{pathwidth}.

For a digraph $G$, we say $P = (W_1,\ldots,W_r)$ is
a {\it path-decomposition} of $G$ if:
\begin{itemize}
\item $r \geq 1~$ and $\cup (W_i:\;1\le i\le r) = V(G)$,
\item {\it (betweenness condition)} for $1\leq h < i < j \leq r$, $W_h \cap W_j \subseteq W_i$, and
\item {\it (cut condition)} if $uv \in E(G)$, then there exist $i,j$ with $1\le i\le j\le r$ such that $u \in W_j$ and $v \in W_i$.
\end{itemize}
The betweenness condition implies that  $\{i\;: v \in W_i \}$ is an integer interval for each $v \in V(G)$, and the cut condition implies that
for $1\le i\le r$ there are no edges from
$\cup_{h<i}W_h$ to $\cup_{j>i}W_j$ in $G \setminus W_i$.

We define the {\it width} of a path-decomposition $P$ as $\max_{1\leq i \leq r}|W_i| - 1$ and denote it by $pw(P)$.
We say $G$ has {\it path-width} at most $k$ if there exists some path-decomposition $P$ of $G$ with $pw(P) \leq k$
and we denote the path-width of $G$ by $pw(G)$.
For example, a non-null loopless digraph $G$ is acyclic if and only if $pw(G)=0$.

For a path-decomposition $P = (W_1, \ldots,W_r)$, we denote
$\min_{1\le i\le r}|W_i|$,
$\max_{1\le i\le r}|W_i|$,
$W_1$, and $W_r$ by
$m(P)$, $M(P)$, $F(P)$, and $L(P)$, respectively.\footnote{For convenience, we use $M(P) (= pw(P) +1)$ instead of $pw(P)$ in later sections.}
We first prove that having bounded path-width is a minor-closed property.

\begin{thm}\label{closed}
If a digraph has path-width at most $k$, then so do all its minors.
\end{thm}

\Proof
Let $P=(W_1,\ldots,W_r)$ be a path-decomposition of a digraph $G$. Then
$P$ is a path-decomposition of $G \setminus e$ for each edge $e \in E(G)$, and
$(W_1 \setminus v, \ldots, W_r \setminus v)$ is a path-decomposition of $G \setminus v$ for each vertex $v \in V(G)$.
Therefore the path-width of a digraph $G$ does not increase by deleting an edge or a vertex.

Thus, it remains to show that if $G$ has path-width at most $k$,
then so does $G/H$ where $H$ is a strongly-connected subdigraph of $G$.
($G/H$ is the digraph obtained from $G$ by contracting $H$ to a single vertex $w$).
Let $P=(W_1,...,W_r)$ be a path-decomposition of $G$, and let $I_H =  \{i\;: W_i \cap V(H) \neq \emptyset\}$.
\\
\\
(1) {\it $I_H$ is an integer interval.}
\\
\\
Suppose $I_H$ is not an integer interval. Take indices $h<i<j$ such that $h,j \in I_H$ and $i \notin I_H$.
Let $u\in W_h \cap V(H)$ and $v\in W_j \cap V(H)$. Since  $\{t\;: u \in W_t \} \subseteq \{1,\ldots, i-1\}$ and $\{t\;: v \in W_t\} \subseteq \{i+1,\ldots,r\}$,
two sets $\{t\;: u \in W_t \}$ and $\{t\;: v \in W_t \}$ do not intersect. 
Since $H$ is strongly-connected and $V(H)\cap W_i = \emptyset$,
there is a directed path from $u$ to $v$ in $G \setminus W_i$.
However, this contradicts the cut condition since there are no edges from
$\cup_{a<i}W_a$ to $\cup_{b>i}W_b$ in $G \setminus W_i$. This proves (1).

\bigskip

Let $P=(W_1,\ldots,W_r)$ be a path-decomposition of $G$ with $pw(P)\leq k$.
Define $W_i'$ by
$$ W_i' =
 \begin{cases}
 (W_i \setminus V(H)) \cup \{w\} & \text {if } i \in I_H\\
 W_i & \text{otherwise}. \\
 \end{cases}
$$
We claim that $P' = (W_1',\ldots,W_r')$ is a path-decomposition of $G/H$ with $pw(P') \leq k$.
The betweenness condition follows from (1).
For the cut condition,
we only need to consider edges incident with $w$ in $G/H$.
For an edge $uw \in E(G/H)$, consider the corresponding edge $uv \in E(G)$.
By the cut condition for $P$, there exist $i \leq j$ such that $W_j \ni u$ and $W_i \ni v$.
Therefore $W_j' \ni u$ and $W_i' \ni w$.
The same argument applies for edges with tail $w$.
Finally, $pw(P') \leq pw(P) \leq k$ from the definition of $P'$.
This proves \ref{closed}.~\bbox

We introduce a notion and a theorem from~\cite{pathwidth}.
For a digraph $G$,
let $A,B$ and $C$
be mutually disjoint subsets of $V(G)$.
We say $(A,B,C)$ is a {\it $k$-triple} if
\begin{itemize}
\item $|A| = |B| = |C| = k$,
\item $ab \in E(G)$ for every $a \in A$ and $b \in B$,
\item $bc \in E(G)$ for every $b \in B$ and $c \in C$, and
\item $A,C$ can be numbered as $\{a_1,\ldots,a_k\}$ and $\{c_1,\ldots,c_k\}$ respectively such that
$c_ia_i \in E(G)$ for $i = 1,\ldots, k$.
\end{itemize}

\begin{thm}\label{triple}
 Let $(A,B,C)$ be a $k$-triple of a digraph $G$. Then $G$ contains every semi-complete digraph with $k$ vertices as a minor.
\end{thm}
\Proof Let $\{a_1,\ldots,a_k\}$ and $\{c_1,\ldots,c_k\}$ be numberings of $A$ and $C$ such that $c_ia_i \in E(G)$ for $i = 1,\ldots,k$.
Take an ordering $\{b_1,\ldots,b_k\}$ of $B$. Then $G|\{a_i, b_i, c_i\}$ is strongly-connected for each $i$.
Let $G'$ be the digraph obtained from $G|(A \cup B \cup C)$ by contracting $G|\{a_i,b_i,c_i\}$ to a single vertex for each $i$.
Then $|V(G')| = k$ and $uv \in E(G')$ for every distinct $u,v \in V(G')$.
Therefore every semi-complete digraph with $k$ vertices is a subgraph of $G'$ and hence, a minor of $G$.
This proves \ref{triple}.~\bbox

The following theorem says a semi-complete digraph has large path-width if and only if it has a large $k$-triple.

\begin{thm}\label{path1}
For every set $\mathcal{S}$ of semi-complete digraphs, the following are equivalent:
\begin{enumerate}
\item There exists $k$ such that every member of $\mathcal{S}$ has path-width at most $k$.
\item There is a digraph $H$ such that no subdivision of $H$ is a subdigraph of any member of $\mathcal{S}$.
\item There exists $k$ such that for each $G \in \mathcal{S}$, there is no $k$-triple in $G$.
\item There exists $k$ such that for each $G \in \mathcal{S}$, there do not exist $k$ vertices of $G$ that are pairwise $k$-connected.
\item There is a digraph $H$ such that for each $G \in \mathcal{S}$, $G$ does not contain $H$ as a minor.
\end{enumerate}
\end{thm}
\Proof The equivalence of the first four statements was proved by Alexandra Fradkin and the second author in~\cite{pathwidth}.
Here, we prove 1 $\Rightarrow$ 5 $\Rightarrow$ 3 to extend the theorem.

Suppose 1 holds for $k$ and $\mathcal{S}$ and let $H$ be a digraph with $pw(H) > k$. Then 5 holds by \ref{closed}.
Now, suppose 5 holds for $H$ and $\mathcal{S}$.
Let $H'$ be a simple digraph containing $H$ as a minor.
Then 5 holds for $H'$ and $\mathcal{S}$ as well.
By \ref{triple}, $G$ has no $|V(H')|$-triple for each $G \in \mathcal{S}$. Therefore 3 holds. This proves \ref{path1}.~\bbox

Thanks to \ref{path1}, it is enough to show the following statement to prove \ref{main1}.

\begin{thm}\label{main2}
For all $k\ge 0$, the class of all semi-complete digraphs with path-width $\leq k$ is a wqo under minor containment.
\end{thm}

\noindent{\bf Proof of \ref{main1}, assuming \ref{main2}.\ \ }
Let $G_1, G_2, \ldots$ be an infinite sequence of semi-complete digraphs.
We may assume $G_i$ does not contain $G_1$ as a minor for each $i \ge 2$.
From \ref{path1}, there exists $k$ such that every member of $\mathcal{S} = \{G_2,G_3,\ldots\}$ has path-width at most $k$.
From \ref{main2}, there exist $j>i\geq 2$ such that $G_i$ is a minor of $G_j$. This proves \ref{main1}.~\bbox

Most of the remainder of this paper is devoted to proving \ref{main2}.
In section 2, we show the existence of a ``linked'' path-decomposition; and
use it in sections 3 and 4 to prove a slightly more general version of $\ref{main2}$.
In section 5, we give counter-examples to the analogue of \ref{main1} for some super-classes of the class of all semi-complete digraphs.

\section{Linked path-decompositions}

In this section, we make a particular choice of path-decomposition,
which will be used in subsequent sections.
Roughly speaking,
we will break this path-decomposition into a sequence of small path-decompositions in the natural way,
so that we can apply Higman's sequence theorem to this sequence.

We first give some definitions. A directed path $P$ in a semi-complete digraph $G$
is {\it induced} if
$v_iv_j \notin E(G)$ for $j-i \geq 2$ where $v_1,\ldots,v_n$ are the vertices of $P$ in order.
Note that $G|V(P)$ is strongly-connected unless it is a one-edge directed path.
For two sets $A$ and $B$, denote by {\it $A \Delta B$} the symmetric difference
$(A \setminus B) \cup (B \setminus A)$.
For a digraph $G$, we say $(C,D)$ is a {\it separation of $G$ of order $s$} if:
\begin{itemize}
\item $C \cup D = V(G)$,
\item $|C \cap D| = s$, and
\item there are no edges $uv \in E(G)$ such that $u \in C \setminus D$ and $v \in D \setminus C$.
\end{itemize}
If $A,B\subseteq V(G)$, a separation $(C,D)$ of $G$ {\em separates} $A,B$ if
$A\subseteq C$ and $B \subseteq D$.
A path-decomposition $P = (W_1,\ldots,W_r)$ of a digraph $G$ is called a
{\it linked path-decomposition} if:
\begin{itemize}
\item {\it (increment condition)} $|W_i ~ \Delta ~ W_{i+1}| = 1 $ for $i = 1,\ldots,r-1$,
\item {\it (cardinality condition)} $|W_1| = |W_r| = m(P)$, and
\item {\it (linked condition)} if $|W_i| \geq t$ for every $i$ with $h \leq i \leq j$, then there exist $t$ vertex-disjoint directed paths from $W_h$ to $W_j$.
\end{itemize}
Observe that for every $v \notin W_1$, there exists a unique $i$ with
$1\leq i \leq r-1$ such that $\{v\} = W_{i+1} \setminus W_i$,
and for every $v \notin W_r$, there exists a unique $j$ with
$1 \leq j \leq r-1$ such that $\{v\} = W_j \setminus W_{j+1}$.
Therefore the increment condition implies that $r -1 = 2 \big(|V(G)| - m(P)\big)$.
In particular, $r$ is bounded above by $2|V(G)|+1$.
We now prove the existence of a linked path-decomposition in a semi-complete digraph.

\begin{thm}\label{linked}
Let $G$ be a semi-complete digraph and $A, B \subseteq V(G)$ with $|A| = |B| = m \geq 0$.
Suppose there exist $m$ vertex-disjoint directed paths from $A$ to $B$ in $G$, and
there exists a path-decomposition (not necessarily linked) $P$ of $G$ such that
$F(P) = A$, $L(P) = B$, and $M(P) \leq k$ for some $k$.
Then there exists a linked path-decomposition $P'$
such that $F(P') = A$, $L(P') = B$, and $M(P') \leq k$.

In particular, every semi-complete digraph $G$ with $pw(G) \leq k$
has a linked path-decomposition $P$ with $pw(P) \leq k$ and $m(P) = 0$.
\end{thm}
\Proof
Observe that
we can obtain a path-decomposition $P'$ with $F(P') = A$, $L(P')=B$ and $M(P') \leq k$ satisfying the increment condition,
by modifying $P$ as follows: we remove (one of) any two consecutive sets that are equal, and insert appropriate sets between sets that differ by more than one vertex.
The cardinality condition for $P'$ is guaranteed by the existence of $m$ vertex-disjoint paths from $A$ to $B$.

Now, among all the path-decompositions $P=(W_1,\ldots,W_r)$ of $G$ with $F(P) = A$, $L(P) = B$, and $M(P) \leq k$ satisfying the increment condition and the cardinality condition,
we take one with $(n_0, n_1, \ldots, n_k)$ ``lexicographically maximal'', where
$n_j = |\{i : |W_i| = j \}|$ for $j= 0,\ldots,k$.
(More precisely, take one with $n_0$ as large as possible; subject to that, take $n_1$ as large as possible; and so on.)
We can take such a path-decomposition since $r$ is bounded above by $2|V(G)|+1$.
Let $P'$ be the path-decomposition we choose. We show $P'$ satisfies the linked condition.

Let $P'= (W_1, \ldots, W_r)$.
Suppose that $|W_i| \geq t$ for every $i$ with $h\leq i \leq j$ and there do not exist $t$ vertex-disjoint directed paths from $W_h$ to $W_j$.
Then from Menger's theorem, there is a separation $(C,D)$ of order less than $t$ that separates $\cup_{i \leq h}W_i$, $\cup_{i \geq j}W_i$.
Take such a separation $(C,D)$ with minimum order $s$.
We claim that there exist two path-decompositions
$$P^C = (W_1^C, \ldots, W_j^C)$$
of $G|C$ with  $W_i^C = W_i$ for $1 \leq i \leq h$, $W_j^C = C \cap D$, and  $M(P^C) \leq k$, and
$$P^D = (W_h^D, \ldots,W_r^D)$$
of $G|D$ with $W_h^D = C \cap D$, $W_i^D = W_i$ for $j \leq i \leq r$, and $M(P^D) \leq k$.
Then we will show that the ``concatenation'' of the two path-decompositions yields a path-decomposition lexicographically better than $P'$, which contradicts our choice of $P'$.

We construct $P^C$ as follows.
Note that there exist $s$ vertex-disjoint paths from $W_h$ to $W_j$ by the minimality of $s$.
Take $s$ vertex-disjoint directed paths $P_1, \ldots, P_s$ from $W_1\cup\cdots\cup W_h$ to $W_j\cup\cdots\cup W_r$ with minimal union.
For $1\le l\le s$, the minimality of the union of $P_1, \ldots, P_s$ implies that
$P_l$ is induced and no vertex of $P_l$ belongs to  $W_1\cup\cdots\cup W_h$ except its first vertex. Since there is no edge from 
$W_1\cup\cdots\cup W_{h-1}$ to $W_{h+1}\cup\cdots\cup W_r$ in $G\setminus W_h$, it follows that the first vertex of $P_l$ belongs to $W_h$.
Similarly, the last vertex of $P_l$ belongs to $W_j$, and no other vertex of $P_l$ belongs to $W_j\cup\cdots\cup W_r$. Let $p_l \in V(P_l) \cap (C \cap D)$.
\\
\\
(1) {\em For $1\le l\le s$, $\{i :  W_i \cap (D \cap V(P_l))\neq \emptyset,~ 1 \leq i \leq j \}$ is an integer interval containing~$j$.}
\\
\\
If $G| (D \cap V(P_l))$ is strongly-connected, then (1) holds by the same argument as (1) in the proof of \ref{closed}.
Therefore we may assume $G| (D \cap V(P_l))$ has exactly two vertices $u,v$ with one edge $uv$, since $P_l$ is induced.
By the cut condition for the edge $uv$, there exist $a, b$ with $1\le a\le b\le r$ such that $u\in W_b$ and $v\in W_a$.
Since no vertex of $P_l$ except its last is in  $W_j\cup\cdots\cup W_r$, it follows that $b<j$.
Since $v\in W_a\cap W_j$ and $a\le b\le j$, it follows that
$v\in W_b$. In summary,
there exists $b$ with $1\le b<j$ such that $u,v\in W_b$.

On the other hand, $\{i : W_i \cap \{u\} \neq \emptyset,~ 1\leq i \leq j\}$ and
$\{i : W_i \cap \{v\} \neq \emptyset,~ 1\leq i \leq j\}$ are both integer intervals by the betweenness condition. Since they intersect,
the set in question is also an integer interval since it is the union of the two intersecting intervals.
This proves (1).

\bigskip

For each $i$ with $1 \leq i \leq j$, define 
$$W_i^C = (W_i \cap C) ~ \cup ~ \{p_l : W_i \cap (D \cap V(P_l)) \neq \emptyset,~ 1\leq l \leq s\}.$$
Let $P^C = (W_1^C, \ldots, W_j^C)$.
\\
\\
(2) {\em $P^C$ is a path-decomposition of $G|C$ with  $W_i^C = W_i$ for $1 \leq i \leq h$, $W_j^C = C \cap D$, and  $M(P^C) \leq k$.}
\\
\\
It is easy to check that $\cup_{i=1}^{r}W_i^C = C$, and
the betweenness condition follows from (1).
For the cut condition, we only need to consider edges incident with $p_l$ in $G|C$ and this is trivial since
$$\{i : p_l \in W_i, 1 \leq i \leq j \}   \subseteq  \{i : p_l \in W_i^C, 1 \leq i \leq j\}.$$
Therefore $P^C$ is a path-decomposition of $G|C$.

For $1 \leq i \leq h$,
$W_i^C = W_i$  since $W_i \subseteq C$.
And $W_j^C = C\cap D$ since $W_j \subseteq D$.
Finally, $M(P_C) \leq M(P) \leq k$ since $|W_i^C| \leq |W_i|$ for every $i$ with $1 \leq i \leq j$. This proves (2).

\bigskip

Similarly, let $P^D = (W_h^D, \ldots, W_r^D)$, where
$$W_i^D = (W_i \cap D) ~ \cup ~ \{p_l : W_i \cap C \cap V(P_l) \neq \emptyset,~ 1\leq l \leq s\}$$
for each $i$ with $h \leq i \leq r$,
then it is a path-decomposition of $G|D$ with $W_h^D = C \cap D$, $W_i^D = W_i$ for $ j \leq i \leq r$, and $M(P_D) \leq k$.

Let $P^*$ be the path-decomposition of $G$ obtained by concatenating $P^C$ and $P^D$ and refining it to satisfy the increment condition.
Then $P^*$ is ``lexicographically better'' than $P'$ since every $W_a$ with $|W_a|\le s$ is a term in the sequence $P^*$ 
(because $W_i^C = W_i$ for $1 \leq i \leq h$ and $W_i^D = W_i$ for $j \leq i \leq r$, and $|W_i|>s$ for $h\le i\le j$), and
there exists at least one more set of size $s$, namely $C \cap D$.
This proves \ref{linked}.~\bbox

\section{Labeled minors}

In this section, for a wqo $Q$ and a semi-complete digraph $G$,
we assign an element of $E(Q)$ to each vertex of $G$, and we fix a linked path-decomposition $P$ of $G$ together with $m(P)$ vertex-disjoint directed paths from $F(P)$ to $L(P)$.
We define a minor relation for these slightly more general objects and prove a well-quasi-order theorem for them. Then \ref{main2} will follow as an corollary.
Roughly speaking, we need this ``$Q$-labeling'' in order to handle the case when one of the induced directed paths has length one so that we may not contract it.

For integers $m,k$ with $k \geq m \geq 0$ and a well-quasi-order $Q$, we say $D = (G,P,R,l)$ is a {\it $(Q,m,k)$-digraph} if:
\begin{itemize}
 \item $G$ is a semi-complete digraph,
 \item $P$ is a linked path-decomposition of $G$ with $m(P) = m$ and $M(P) \leq k$,
 \item $R = (R_1,\ldots,R_m)$ is a sequence of $m$ vertex-disjoint induced directed paths from $F(P)$ to $L(P)$ in $G$, and
 \item $l$ is a mapping from $V(G)$ to $E(Q)$.
\end{itemize}
Note that $|F(P) \cap V(R_i)| = |L(P) \cap V(R_i)| = 1$ for each $i = 1,\ldots, m$.
We say the vertex in $F(P) \cap V(R_i)$ is the {\it $i$-th source root} of $D$
and the vertex in $L(P) \cap V(R_i)$ is the {\it $i$-th terminal root} of $D$.
We denote the collection of all $(Q,m,k)$-digraphs by $\mathcal{G}_m^k(Q)$.
We say $D$ is {\it trivial} if $r=1$ where $P=(W_1,\ldots,W_r)$.
Note that $|V(G)| = m$ if $D$ is trivial.
\bigskip

Now we define a minor relation on $\mathcal{G}_m^k(Q)$.
Let $D = (G,P,R,l)$ , $D' = (G',P',R',l')$ $\in \mathcal{G}_m^k(Q)$.
Let $a_i$ and $a_i'$ be the $i$-th source roots of $D$ and $D'$, respectively, and similarly let $b_i$ and $b_i'$ be the $i$-th terminal roots of $D$ and $D'$, respectively.
We say $D$ is a {\it minor} of $D'$ if there exists a minor mapping $\phi$ from $G$ to $G'$ such that:
\begin{itemize}
 \item $a_i' \in V(\phi(a_i))$ and  $b_i' \in V(\phi(b_i))$ for $i=1,\ldots,m$, and
 \item for every $v \in V(G)$, $l(v) \leq_Q l'(u)$ for some $u \in V(\phi(v))$.
\end{itemize}
Again, we call $\phi$ a {\it minor mapping} from $D$ to $D'$.

Next, we define a ``decomposition'' of a $(Q,m,k)$-digraph.
Let $D=(G,P,R,l)$ be a $(Q,m,k)$-digraph with $P = (W_1,\ldots,W_r)$, $R=(R_1,\ldots,R_m)$ and suppose $|W_s| = m$ for some $s$ with $1 < s < r$.
Let $A = \cup_{i\leq s}W_i$ and define $D|A = (G_A,P_A,R_A,l_A)$ by:
\begin{itemize}
 \item $G_A = G|A$,
 \item $P_A = (W_1,\ldots,W_s)$,
 \item $R_A = (R_1 |A, \ldots, R_m |A)$, and
 \item $l_A = l | A$.
\end{itemize}
Then $D|A$ is a $(Q,m,k)$-digraph and similarly, $D|B$ is a $(Q,m,k)$-digraph where $B = \cup_{i \geq s}W_i$.
We write $D = D_A \oplus D_B$ and say $D$ is {\it decomposable}.
We denote the class of all non-trivial non-decomposable $(Q,m,k)$-digraphs by $\mathcal{ND}_m^k(Q)$.
More precisely, we say $D=(G,P,R,l) \in \mathcal{G}_m^k(Q)$ is in $\mathcal{ND}_m^k(Q)$ if:
\begin{itemize}
 \item $r \geq 3$, and $|W_i| > m$ for every $i$ with $2 \leq i \leq r-1$ where $P = (W_1,\ldots,W_r)$.
\end{itemize}
Note that $P' = (W_2,\ldots,W_{r-1})$ is a linked path-decomposition of $G$ with $m(P') = m+1$, $F(P') = W_2 \supseteq W_1$, and $L(P') = W_{r-1} \supseteq W_r$.
Let $R'$ be a sequence of $m+1$ vertex-disjoint induced directed paths from $W_2$ to $W_{r-1}$.
Then we see that each $D=(G,P,R,l) \in \mathcal{ND}_m^k(Q)$ yields at least one member $D'=(G,P',R',l) \in \mathcal{G}_{m+1}^k$.
(Notice that it could be the case that some path in $R'$ joins the $i$-th source root of $D$ to
the $j$-th terminal root of $D$ for some $j \neq i$.)

\begin{thm}\label{extendable}
Let $m,k$ be integers with $k > m \geq 0$. Suppose $\mathcal{G}_{m+1}^{k}(Q)$ is a wqo under minor containment for every wqo $Q$.
Then $\mathcal{ND}_m^k(Q)$ is a wqo under minor containment for every wqo $Q$ as well.
\end{thm}
\Proof
Let $Q$ be a wqo and $D_1,D_2,\ldots$ be an infinite sequence in $\mathcal{ND}_m^k(Q)$.
For each $D_i = (G_i,P_i,R_i,l_i)$, let $D_i'=(G_i,P_i',R_i',l_i) \in \mathcal{G}_{m+1}^{k}(Q)$ as described earlier.
Recall that every source root of $D_i$ is also a source root of $D_i'$, and every terminal root of $D_i$ is also a terminal root of $D_i'$.

For each $i \geq 1$, let $\sigma_i , \tau_i: \{1,\ldots,m\} \rightarrow \{1,\ldots,m+1\}$ be injections defined by
\begin{itemize}
 \item the $t$-th source root of $D_i$ equals the $\sigma_i(t)$-th source root of $D_i'$.
 \item the $t$-th terminal root of $D_i$ equals the $\tau_i(t)$-th terminal root of $D_i'$.
\end{itemize}
Since there are only finitely many pairs $(\sigma_i,\tau_i)$,
there exists some $(\sigma,\tau)$ such that $(\sigma_i,\tau_i) = (\sigma,\tau)$ for infinitely many $i$.
Therefore we may assume $\sigma_i = \sigma$ and $\tau_i = \tau$ for every $i\geq 1$.

Since $\mathcal{G}_{m+1}^{k}(Q)$ is a wqo under minor containment,
$D_i'$ is a minor of $D_j'$  for some $i<j$ with some minor mapping $\phi$.
Then $\phi$ is also a minor mapping from $D_i$ to $D_j$.
This proves \ref{extendable}.~\bbox

We say $D = (G,P,R,l) \in \mathcal{G}_m^k(Q)$ is {\it contractible} if
\begin{itemize}
 \item $G|V(R_j)$ is strongly-connected for every $j \in \{1,\ldots,m\}$ where $R = (R_1,\ldots,R_m)$.
\end{itemize}
We denote the set of all non-contractible $(Q,m,k)$-digraphs by $\mathcal{NC}_m^k(Q)$.

Suppose $G \in \mathcal{NC}_m^k(Q)$. In other words, $G|V(R_j)$ is not strongly-connected for some $j$.
Then $G|V(R_j)$ must be a digraph with two vertices, namely the $j$-th source root $u$
and the $j$-th terminal root $v$, and one edge $uv$.
Note that every $W_i$ contains either $u$ or $v$ where $P = (W_1,\ldots,W_r)$.
Let $G' = G \setminus \{u,v\}$,
$\hat{W_i} = W_i \setminus \{u,v\}$ for $i=1,\ldots,r$, and $\hat{P} = (\hat{W_1},\ldots,\hat{W_r})$.
Then $\hat{P}$ is a path-decomposition (not necessarily linked) of $G'$ with  $M(\hat{P}) \leq k-1$.
Note that we still have $m-1$ vertex-disjoint paths from $F(P) \setminus u$ to $L(P) \setminus v$.
From \ref{linked}, there exists a linked path-decomposition $P'$ of $G'$ with $F(P') = F(P) \setminus u$, $L(P') = L(P) \setminus v$,
$m(P') = m-1$, and $M(P') \leq k-1$.
Also, the sequence $R'$ obtained from $R$ by omitting $R_j$ is a sequence of $m-1$ vertex-disjoint induced directed paths from $F(P')$ to $L(P')$.
For labels, let $Q'$ be a well-quasi-order defined by
\begin{itemize}
\item $E(Q') = E(Q) \times \{0,1,2\} \times \{0,1,2\}$, and
\item $(q,x,y) \leq_{Q'} (q',x',y')$ if and only if $q \leq_Q q'$ and $x = x'$ and $y = y'$.
\end{itemize}
For each $w \in V(G) \setminus \{u,v\}$, let

$$ x(w) =
 \begin{cases}
 0 & wu \in E(G), uw \notin E(G)\\
 1 & wu \notin E(G), uw \in E(G)\\
 2 & wu \in E(G), uw \in E(G)\\
 \end{cases}
$$

$$ y(w) =
 \begin{cases}
 0 & wv \in E(G), vw \notin E(G)\\
 1 & wv \notin E(G), vw \in E(G)\\
 2 & wv \in E(G), vw \in E(G).\\
 \end{cases}
$$
Let $l'$ be a mapping from $V(G')$ to  $E(Q')$ defined by
\begin{itemize}
 \item $l'(w) = (l(w),x(w),y(w))$ for each $w \in V(G) \setminus \{u,v\}$.
\end{itemize}
Then $D' = (G',P',R',l') \in \mathcal{G}_{m-1}^{k-1}(Q')$ and
we see that each $D \in \mathcal{ND}_m^k(Q)$ yields at least one member $D'$ in $\mathcal{G}_{m-1}^{k-1}(Q')$.

\bigskip

\begin{thm}\label{almost}
Let $m,k$ be integers with $k \geq m \geq 1$. Suppose $\mathcal{G}_{m-1}^{k-1}(Q)$ is a wqo under minor containment for every wqo $Q$.
Then $\mathcal{NC}_m^k(Q)$ is a wqo under minor containment for every wqo $Q$ as well.
\end{thm}
\Proof
Let $Q$ be a wqo and $D_1,D_2,\ldots$ be an infinite sequence in $\mathcal{NC}_m^k(Q)$.
For each $D_i = (G_i,P_i,R_i,l_i)$, let $D_i' = (G_i',P_i',R_i',l_i') \in \mathcal{G}_{m-1}^{k-1}(Q')$ as described earlier
and let $u_i$ and  $v_i$ be the source root and the terminal root of $(G_i,P_i,R_i,l_i)$ such that $G_i' = G_i \setminus \{u_i,v_i\}$.
Since $Q$ is a wqo, we may assume
\begin{itemize}
 \item $l_i(u_i) \leq_Q l_j(u_j)$, and $l_i(v_i) \leq_Q l_j(v_j)$ for every $i<j$.
\end{itemize}
Since $\mathcal{G}_{m-1}^{k-1}(Q')$ is a wqo under minor containment, there exist $i,j$ with $1\leq i < j$ such that
$D_i'$ is a minor of $D_j'$  with a minor mapping $\phi'$.
Define $\phi$ from $D_i$ to $D_j$ as
$$ \phi(w) =
 \begin{cases}
 (\{u_j\},\emptyset) & w=u_i \\
 (\{v_j\},\emptyset) & w=v_i \\
 \phi'(w) & w \in V(G_i')\\
 \end{cases}
$$
Then it is easy to check that $\phi$ is a minor mapping from $D_i$ to $D_j$, by the definition of $Q'$-labels.
This proves \ref{almost}.~\bbox

For two subclasses $\mathcal{A}, \mathcal{B}$ of $\mathcal{G}_m^k(Q)$,
denote by $\mathcal{A} \oplus \mathcal{B}$ the class of all $(Q,m,k)$-digraphs $D$ which are decomposable as $D_A \oplus D_B$
where $D_A \in \mathcal{A}$ and $D_B \in \mathcal{B}$.

\begin{thm}\label{union}
If $\mathcal{A} , \mathcal{B} \subseteq \mathcal{G}_m^k(Q)$ are both wqo under minor containment, then $\mathcal{A} \cup \mathcal{B}$ and $\mathcal{A} \oplus \mathcal{B}$ are both wqo under minor containment.
\end{thm}
\Proof
$\mathcal{A} \cup \mathcal{B}$ is a wqo under minor containment because every infinite sequence in $\mathcal{A} \cup \mathcal{B}$
contains either an infinite subsequence in $\mathcal{A}$ or an infinite subsequence in $\mathcal{B}$.

For $\mathcal{A} \oplus \mathcal{B}$,
let $D_1,D_2,\ldots$ be an infinite sequence in $\mathcal{A} \oplus \mathcal{B}$. Let
$D_i = D_i^a \oplus D_i^b$
where $D_i^a = (G_i^a,P_i^a,R_i^a,l_i^a) \in \mathcal{A}$ and $D_i^b = (G_i^b,P_i^b,R_i^b,l_i^b) \in \mathcal{B}$ for each $i \geq 1$.
Then there exist $i<j$ such that:
\begin{itemize}
 \item $D_i^a$ is a minor of $D_j^a$ with minor mapping $\phi_a$, and
 \item $D_i^b$ is a minor of $D_j^b$ with minor mapping $\phi_b$.
\end{itemize}
Define a mapping $\phi$ from $D_i$ to $D_j$ by
$$ \phi(w) =
 \begin{cases}
 \phi_a(w) & w \in V(G_i^a) \setminus V(G_i^b) \\
 \phi_b(w) & w \in V(G_i^b) \setminus V(G_i^a) \\
 \phi_a(w) \cup \phi_b(w) & w \in V(G_i^a) \cap V(G_i^b).\\
 \end{cases}
$$
Since the union of two strongly-connected subdigraphs with non-empty intersection is also strongly-connected, $\phi_a(w) \cup \phi_b(w)$ is strongly-connected in $G_j$ for
each $w \in V(G_i^a) \cap V(G_i^b)$. Then it is easy to check that $\phi$ is a minor mapping from $D_i$ to $D_j$.
This proves \ref{union}.~\bbox

\section{Links, and the main proof}

We say a contractible $(Q,m,k)$-digraph $D$ ($\notin \mathcal{NC}_m^k(Q)$) is a {\it link} if:
\begin{itemize}
 \item $D \in  \mathcal{ND}_m^k(Q) ~\cup ~(\mathcal{NC}_m^k(Q) \oplus \mathcal{ND}_m^{k}(Q)) $.
\end{itemize}
We denote the collection of all links in $\mathcal{G}_m^k(Q)$ by $\mathcal{L}_m^k(Q)$.
The following is an easy corollary of \ref{union}.

\begin{thm}\label{link}
Suppose $\mathcal{NC}_m^k(Q)$ and $\mathcal{ND}_m^{k}(Q)$ are wqo under minor containment for some wqo $Q$.
Then $\mathcal{L}_m^k(Q)$ is a wqo under minor containment as well.
\end{thm}

Now, we decompose $D \in \mathcal{G}_m^k(Q)$ into links (possibly except the last term) to apply Higman's sequence theorem.

\begin{thm}\label{concatenation}
Let $D=(G,P,R,l)$ be a non-trivial $(Q,m,k)$-digraph. Then
$D = D_1 \oplus \ldots \oplus D_t$ (perhaps with $t = 1$) such that:
\begin{itemize}
 \item $D_i \in \mathcal{L}_m^k(Q)$ for $i \in \{1,\ldots,t-1\}$, and
 \item $D_t \in \mathcal{L}_m^k(Q)~ \cup ~\mathcal{NC}_m^k(Q)$.
\end{itemize}
\end{thm}
\Proof
We may assume $D$ is contractible since otherwise $D\in \mathcal{NC}_m^k(Q)$ and the result holds with $t = 1$.
Let $P = (W_1,\ldots,W_r)$. We proceed by induction on $r$.
For the base case $r=3$, $D$ belongs to $\mathcal{ND}_m^k(Q)$ and hence, $D$ itself is a link.

Let $1 = n_1 < \cdots < n_s = r$ be the indices such that $$|W_{n_1}| = \cdots = |W_{n_s}| = m$$
Let $j > 1$ be the smallest index such that the initial segment $D | (\cup_{i=1}^{n_j}W_i)$ is contractible
(such $j$ exists since $D$ is contractible).

If $j = 2$, then the initial segment $D | (\cup_{i=1}^{n_2}W_i)$ belongs to $\mathcal{ND}_m^k(Q)$, and hence it is a link.
If $j > 2$, then
$$D | (\cup_{i=1}^{n_j}W_i)
= D | (\cup_{i=1}^{n_{j-1}}W_i)~ \oplus ~ D | (\cup_{i=n_{j-1}}^{n_{j}}W_i)$$
$$\in \mathcal{NC}_m^k(Q)~ \oplus ~\mathcal{ND}_{m}^{k}(Q).$$
Therefore, in either case, $D | (\cup_{i=1}^{n_j}W_i)$ is a link.
If $j = s$, then $D$ is a link, and we are done. Otherwise,
$D | (\cup_{i=n_j}^{r}W_i)$ is non-trivial and satisfies the statement by the induction hypothesis.
Therefore $$D = D | (\cup_{i=1}^{n_j}W_i) ~ \oplus ~ D | (\cup_{i=n_j}^{r}W_i).$$
satisfies the statement as well.
This proves \ref{concatenation}.~\bbox

Let $Q$ be a quasi-order. We define a quasi-order $Q^{<\omega}$ on the set of all finite sequences of elements of $E(Q)$.
Let $p = (p_1,\ldots,p_a)$ and $q = (q_1,\ldots,q_b)$ be sequences of elements of $E(Q)$.
Then $p \leq_{Q^{<\omega}} q$ if and only if:
\begin{itemize}
 \item $a \leq b$, and
 \item there exist $1 \leq \alpha_1 <\ldots<\alpha_a \leq b$ such that $p_i \leq_Q q_{\alpha_i}$ for every $i=1,\ldots,a$.
\end{itemize}
It is proved in \cite{Higman} that

\begin{thm}\label{Higman1}
 If $Q$ is a wqo, then so is $Q^{<\omega}$.
\end{thm}

The following is an easy corollary of \ref{Higman1}.

\begin{thm}\label{Higman2}
Suppose $Q_1,Q_2$ and $Q_3$ are wqo. Let $Q$ be a quasi-order with $E(Q)$ the set of all finite sequences $(p_1,\ldots,p_a)$ with $a\ge 2$ such that $p_1\in
E(Q_1)$, $p_2,\ldots,p_{a-1}\in E(Q_2)$ and $p_a\in E(Q_3)$.
For $p = (p_1,\ldots,p_a)$ and $q = (q_1,\ldots,q_b) \in E(Q)$, let
$p \leq_{Q} q$ if and only if:
\begin{itemize}
 \item $a \leq b$, and
 \item there exist $1 = \alpha_1 <\ldots<\alpha_a = b$ such that $p_1 \leq_{Q_1} q_1$, $p_a \leq_{Q_3} q_b$, and $p_i \leq_{Q_2} q_{\alpha_i}$ for every $i=2,\ldots,a-1$.
\end{itemize}
Then $Q$ is a wqo.
\end{thm}

Next, we prove a key lemma for \ref{main2}.

\begin{thm}\label{mainlemma}
Let $m,k$ be integers with $k > m \geq 1$. Suppose $\mathcal{G}_{m-1}^{k-1}(Q)$ and $\mathcal{G}_{m+1}^{k}(Q)$ are both wqo under minor containment for every wqo $Q$.
Then $\mathcal{G}_m^k(Q)$ is a wqo under minor containment for every wqo $Q$.
\end{thm}
\Proof
Let $D_1,D_2,\ldots$ be an infinite sequence of $\mathcal{G}_m^k(Q)$.
We may assume $D_i$ is non-trivial for every $i \geq 1$ because every trivial $D_i$ has $m$ vertices.
Decompose each $D_i$ as $$D_i = D_i^1 \oplus \ldots \oplus D_i^{t_i}$$ as in \ref{concatenation}.
By \ref{extendable}, \ref{almost}, and \ref{link}, $\mathcal{L}_m^k(Q)$ and $\mathcal{NC}_m^k(Q)$ are both well-quasi-ordered under minor containment.
Therefore we may assume $t_i \geq 3$ for every $i \geq 1$.
We apply \ref{Higman2} for $Q_1 = Q_2 = \mathcal{L}_m^k(Q)$, $Q_3 = \mathcal{L}_m^k(Q)~ \cup ~\mathcal{NC}_m^k(Q)$.
Then there exist $i<j$ with $t_i \leq t_j$ and $1=\alpha_1 < \ldots < \alpha_{t_i} =t_j$ such that:
\begin{itemize}
 \item $D_i^p$ is a minor of $D_j^{\alpha_p}$ with a minor mapping $\phi_p$ for every $1 \leq p \leq t_i$ .
\end{itemize}
For each $w \in V(G_i)$, let $w_l \leq w_r$ be the indices such that
$w_l = \min\{l : w \in V(G_i^l)  \}$, and
$w_r = \max\{r : w \in V(G_i^r) \}$. 
If $w_l < w_r$, then $w$ must be on one of the $m$ vertex-disjoint induced directed paths of $D_i$, and let $R_w$ be the path.

Now, define $\phi$ from $D_i$ to $D_j$ by
$$ \phi(w) = \bigcup_{p \in \{w_l,\ldots,w_r\}} \phi_p(w)  \bigcup_{\alpha_{w_l}< \alpha < \alpha_{w_r}} G_j^{\alpha}|V(R_w)$$

Note that $G_j^{\alpha}|V(R_w)$ is strongly-connected since $D_j^{\alpha}$ is a link.
Then $\phi$ is a minor mapping from $D_i$ to $D_j$.
This proves \ref{mainlemma}.~\bbox

\begin{thm}\label{maintheorem}
 Let $m,k$ be integers with $k\geq m\geq 0$. Then $\mathcal{G}_m^k(Q)$ is a wqo under minor containment for every wqo $Q$.
\end{thm}
\Proof
We proceed by induction on $k$. For fixed $k$, we use induction on $k-m$.
For the base case $k = m$, every member $D \in \mathcal{G}_m^k(Q)$ is trivial and hence the statement holds.
For the inductive step, since $\mathcal{G}_{m-1}^{k-1}(Q)$ and $\mathcal{G}_{m+1}^{k}(Q)$ are wqo by the inductive hypotheses, the statement follows from \ref{mainlemma}.
This proves \ref{maintheorem}.~\bbox

\noindent{\bf Proof of \ref{main2}.\ \ }
Let $G_1, G_2, \ldots$ be an infinite sequence of semi-complete digraphs with path-width at most $k$.
Let $Q$ be a wqo with $E(Q) = \{0\}$.
For each $i\ge 1$, 
let $P_i$ be a linked path-decomposition of $G_i$ with $m(P_i) = 0$ and $pw(P_i) \leq k$. 
Let $R_i = ()$ be the empty sequence, let $l_i$ be the constant mapping from $V(G_i)$ to $\{0\}$, and let $D_i = (G_i,P_i,R_i,l_i)$.
Since $D_i \in \mathcal{G}_0^{k+1}(Q)$ for each $i \geq 1$,
there exist $j > i \geq 1$ such that $D_i$ is a minor of $D_j$ by \ref{maintheorem}.
Therefore $G_i$ is a minor of $G_j$.
This proves \ref{main2}.~\bbox

\section{Counter-examples}

In this section, we give some counter-examples for some classes of digraphs containing all semi-complete digraphs.

A digraph $G$ is a {\it super-tournament} if either $uv \in E(G)$ or $vu \in E(G)$ for every distinct $u,v \in V(G)$.
In particular, a simple super-tournament is a semi-complete digraph.
We give a counter-example to show that the class of all super-tournaments is not a wqo under minor containment; and indeed, the subclass of all
super-tournaments with no three edges mutually parallel is not a wqo.

For $i \geq 3$, let $T_i$ be a transitive tournament with $i$ vertices $v_1, \ldots, v_{i}$ such that $v_av_b \in E(T_i)$ if and only if $a < b$. 
Let $G_i$ be a super-tournament obtained from $T_i$ by doubling the following $i$ edges:
$$v_1v_2,~ v_2v_3,~ \ldots, ~v_{i-1},v_i,~ v_1v_i.$$
(``Doubling'' means adding a new edge with the same head and tail as the given edge.)
We claim that $G_i$ is not a minor of $G_j$ for $j>i\geq 3$.
First, we cannot contract anything from $G_j$ because $G_j$ has no directed cycles.
Therefore $G_j$ must have $G_i$ as a subdigraph in order to contain it as a minor.
However, note that the underlying undirected graph of $G_i$ has a cycle of length $i$ with all edges doubled, while
$G_j$ does not.
Therefore $G_i$ is not a subdigraph of $G_j$ and hence, not a minor of $G_j$.

The {\it stability number} for a digraph $G$ is the maximum size of an independent set in the underlying undirected graph of $G$.
For example, a non-null semi-complete digraph has stability number one.
We give a counter-example for the class of simple digraphs with stability number at most two.

For $i \geq 2$,
let $A_i = \{a_1,a_2,a_3\}$, $B_i = \{b_1,b_2,b_3\}$, $C_i = \{c_1,\ldots,c_i\}$, and $D_i = \{d_1,\ldots,d_i\}$.
Let $G_i$ be a simple digraph with stability number two defined as follows. (See figure 1.)
\begin{itemize}
 \item $V(G_i)$ is the disjoint union of $A_i,B_i,C_i$ and $D_i$,
 \item $G_i|A_i$, $G_i|B_i$ are directed triangles,
 \item $G_i|C_i$, $G_i|D_i$ are transitive tournaments,
 \item $A_i$ is complete to $C_i$,
 \item $D_i$ is complete to $B_i$,
 \item $b_1a_1$ is the only edge between $A_i$ and $B_i$.
 \item Each edge between $C_i$ and $D_i$ goes from $C_i$ to $D_i$, and the bipartite graph underlying $(C_i\cup D_i, \delta^+(C_i,D_i))$ is a Hamiltonian cycle.
 \item There are no other edges between $A \cup C_i$ and $B \cup D_i$.
\end{itemize}
We claim that there do not exist $j>i\ge 2$ such that $G_i$ is a minor of $G_j$.
For suppose $G_i$ is a minor of $G_j$ with minor mapping $\phi$.
First, observe the following fact.
\begin{itemize}
 \item If $H$ is a strongly-connected subdigraph of $G_j$ with $|V(H)| \geq 2$, then either $A_j \subseteq V(H)$ or $B_j \subseteq V(H)$ or $b_1a_1 \in E(H)$.
\end{itemize}
Therefore once we contract a non-trivial strongly-connected subdigraph $H$ of $G_j$, there do not exist two disjoint directed cycles. 
That means we cannot contract anything in $G_j$ if we hope to obtain $G_i$ as a minor.
Therefore $\phi$ must be a subdigraph mapping and
$\phi(A_i) = G_j|A_j$ and $\phi(B_i)= G_j|B_j$ in order to preserve the existence of two disjoint directed cycles with an edge between them.
Then $\phi(C_i)$ and $\phi(D_i)$ are subdigraphs of $G_j|C_j$ and $G_j|D_j$, respectively.
However, the underlying bipartite graph of $(C_j\cup D_j, \delta^+(C_j,D_j))$ is a cycle of length $2j$,
and hence does not contain a cycle of length $2i$ as a subgraph, a contradiction.
This proves our claim.

\begin{figure}
  \begin{center}
    \scalebox{0.6}
      {
        \includegraphics{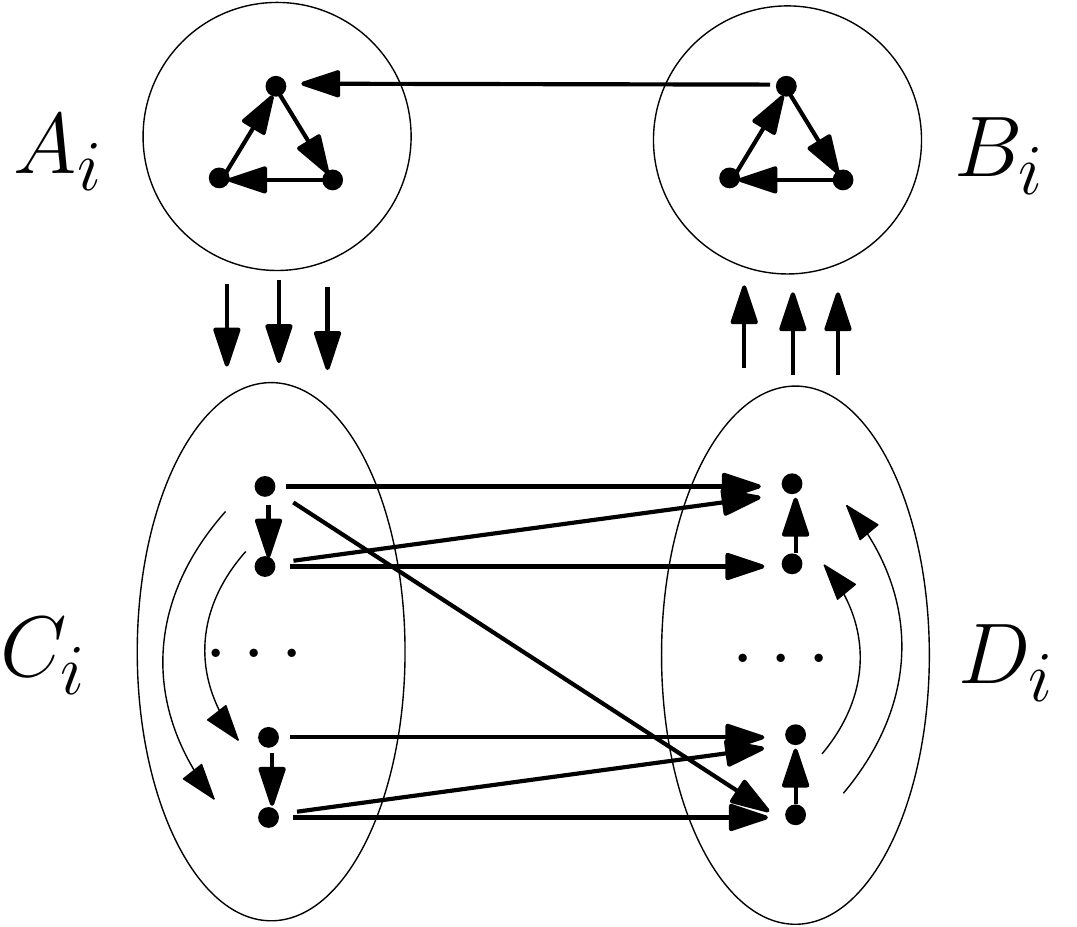}
      }
  \end{center}
  \caption{$G_i$}
\end{figure}

\end{document}